\documentclass{amsart}

\usepackage[all]{xy}
\usepackage{amsmath}
\usepackage{amssymb}
\usepackage{amsthm}
\usepackage{latexsym}
\usepackage{enumerate}
\usepackage{prettyref}
\usepackage{hyperref}
\usepackage{color}

\newtheorem{mainthm}{Theorem}
\newtheorem{conjecture}[mainthm]{Conjecture}
\newrefformat{conj}{\hyperref[{#1}]{Conjecture~\ref*{#1}}}

\newtheorem{theorem}{Theorem}[section]
\newrefformat{thm}{\hyperref[{#1}]{Theorem~\ref*{#1}}}
\newtheorem{definition}[theorem]{Definition}
\newrefformat{def}{\hyperref[{#1}]{Definition~\ref*{#1}}}

\newrefformat{lem}{\hyperref[{#1}]{Lemma~\ref*{#1}}}

\newrefformat{prop}{\hyperref[{#1}]{Proposition~\ref*{#1}}}
\newtheorem{corollary}[theorem]{Corollary}
\newrefformat{cor}{\hyperref[{#1}]{Corollary~\ref*{#1}}}

\newrefformat{rem}{\hyperref[{#1}]{Remark~\ref*{#1}}}
{
\newtheorem{examplecore}[theorem]{Example}}
\newrefformat{ex}{\hyperref[{#1}]{Example~\ref*{#1}}}

\newcommand{\op}{\operatorname}

\begin{document}

\title{Homology of $\op{GL}_3$ of function rings of elliptic curves}

\author{Matthias Wendt}
\address{Matthias Wendt, Fakult\"at Mathematik,
Universit\"at Duisburg-Essen, Thea-Leymann-Strasse 9, 
Essen, Germany} 
\email{matthias.wendt@uni-due.de}

\begin{abstract}
The note provides a description of the homology of $\op{GL}_3$ over
function rings of affine elliptic curves over arbitrary fields,
following the earlier work of Takahashi and Knudson in the case
$\op{GL}_2$. Some prospects for applications to K-theory of elliptic
curves are also discussed.
\end{abstract}

\maketitle

\section{Introduction}

This note is a step towards understanding the homology groups
$\op{H}_\bullet(\op{GL}_n(k[C]),\mathbb{Z})$ of general linear groups
over function rings of affine elliptic curves (with one point $\op{O}$
at infinity). At present, only the
case $\op{(P)GL}_2$ is sufficiently understood, due to the work of
Takahashi \cite{takahashi:tree} and Knudson \cite{knudson:elliptic},
see in particular \cite[Section 4.5]{knudson:book}. Takahashi computed
the action of 
$\op{GL}_2(k[C])$ on the Bruhat-Tits tree associated to $K=k(C)$ with
the valuation corresponding to the point $\op{O}$ at infinity and
found a fundamental domain. Knudson then worked out the equivariant
spectral sequence and computed
$\op{H}_\bullet(\op{PGL}_2(k[C]),\mathbb{Z})$. 

This note takes the next step, working out both the equivariant cell
structure of the $\op{GL}_3(k[C])$-action on the building and
evaluating the resulting equivariant spectral sequence. This can be
done completely explicitly, via matrix calculations as in
\cite{takahashi:tree} (and this actually is the way I first hit on the
results in the fall of 2013). However, in the meantime there is a more
conceptual way of formulating (and proving) the results which should
also be applicable to higher rank computations. 

Let me outline the results of Knudson and Takahashi from this more
conceptual perspective. The quotient
$\op{GL}_2(k[C])\backslash\mathfrak{T}_C$ 
of the Bruhat-Tits tree is related to rank two vector bundles on the
complete elliptic curve $\overline{C}$. The quotient is contractible,
and can be retracted onto (the star of) the central point $o$ of
Takahashi's fundamental domain which corresponds to the unique stable
bundle $E(2,1)$ on $\overline{C}$ which restricts trivially to
$C$. After removing this central point, the fundamental 
domain decomposes into disjoint trees indexed by $\mathbb{P}^1(k)$
corresponding to semistable rank two bundles
$\mathcal{L}\oplus\mathcal{L}^{-1}$ with trivial determinant on
$\overline{C}$. On each of these trees, the homology of the
automorphism groups of the vector bundles is constant (using homotopy
invariance), and we can rewrite Knudson's formula \cite[Theorem
4.5.2]{knudson:book} as
$$
\op{H}_i(\op{PGL}_2(k[C]),\mathbb{Z})\cong
\bigoplus_{E\in\mathcal{M}_{2,\mathcal{O}}(\overline{C})}
\op{H}_i(\op{Aut}(E)/k^\times,\mathbb{Z}). 
$$
The results for $\op{GL}_3(k[C])$ follow the same path: the quotient
$\op{GL}_3(k[C])\backslash\mathfrak{X}_C$ is described in terms of
vector bundles; the structure of the subcomplex of bundles with
non-trivial automorphisms can be understood in terms of suitable
moduli spaces of vector bundles; and the remaining part of the
quotient is controlled by the stable bundles.

There are two important points where the results for $\op{GL}_3$ differ
from those for $\op{GL}_2$. First, the subcomplex of unstable bundles
is essentially zero-dimensional for $\op{GL}_2$, but is essentially
1-dimensional in the case $\op{GL}_3$. Second, the quotient
$\op{GL}_2(k[C])\backslash\mathfrak{T}_C$ is contractible because
there is a unique stable bundle $E(2,1)$, but the quotient
$\op{GL}_3(k[C])\backslash\mathfrak{X}_C$ is not contractible because
there are two stable bundles $E(3,1)$ and $E(3,2)$. 

\emph{Acknowledgement:} The motivation for this work goes back to
early 2011 when Guido Kings asked me if Knudson's work on K-theory of
elliptic curves using buildings could be pushed further. It took me a
while to realize what interesting ramifications this question had, and
it will take me a while longer to work out the answer. The present
note is a snapshot of work in progress...

\section{Buildings and vector bundles}

In the following, let $k$ be a field, $\overline{C}$ be an elliptic
curve with $k$-rational point $\op{O}$, and set
$C=\overline{C}\setminus\{\op{O}\}$. The function field $K=k(C)$ has a
valuation $v=v_{\op{O}}$, and   $\mathfrak{X}_C$ will denote the
Bruhat-Tits building associated to $\op{GL}_3$ and the valuation
$v$. Recall that the building $\mathfrak{X}_C$ is a
contractible simplicial complex whose set of $0$-simplices is the set
$\op{GL}_3(K)/(K^\times\cdot \op{GL}_3(\mathcal{O}_v))$ of equivalence
classes of lattices, and whose higher simplices are given by chains of
lattice inclusions. 

The group $\op{GL}_3(k[C])$ acts on the building $\mathfrak{X}_C$, and
the $0$-simplices of the quotient
$\op{GL}_3(k[C])\backslash\mathfrak{X}_C$ can be identified with the
set of rank three vector bundles on $\overline{C}$ which restrict
trivially to $C$ modulo tensoring with the ideal sheaf
$\mathcal{I}_{\op{O}}$. This identification is well-known, in the case
of $\op{GL}_2$ it can be found in \cite{serre:book}. Under this
identification, the stabilizer of the $0$-simplex corresponding to the
vector bundle $E$ on $\overline{C}$ is the automorphism group
$\op{Aut}(E)$. 

Recall Atiyah's classification \cite{atiyah:classification}
of vector bundles on elliptic curves over algebraically closed fields,
which was extended to arbitrary fields by Pumpl\"un
\cite{pumpluen:classification}. There is a unique indecomposable
bundle $E(r,d)$ of rank $r$ with determinant
$\mathcal{I}_{\op{O}}^{\otimes d}$ (and all indecomposable bundles can
be obtained this way by replacing $\mathcal{I}_{\op{O}}$ by arbitary
line bundles $\mathcal{L}$). These bundles are stable if
$(r,d)=1$, and semi-stable otherwise. Moreover, the Harder-Narasimhan
filtration for vector bundles over elliptic curves splits, so that
each vector bundle is S-equivalent to a direct sum of indecomposable
bundles as above.  

Finally, recall from \cite{tu:semistable} that the moduli space
$\mathcal{M}_{r,d}(\overline{C})$ of S-equivalence classes of
semistable vector bundles of rank $r$ and degree $d$ on $\overline{C}$
is isomorphic to $\op{Sym}^h\overline{C}$ with $h=(r,d)$. Under this
isomorphism, the determinant map
$\det:\mathcal{M}_{r,d}(\overline{C})\to\op{Jac}_d(\overline{C})$ is
identified with the Abel-Jacobi map, so that for a fixed line bundle
$\mathcal{L}$, the moduli space $\mathcal{M}_{r,\mathcal{L}}$ of
semistable rank $r$ vector bundles with determinant $\mathcal{L}$ is
identified with $\mathbb{P}^{h-1}$. 

\section{Unstable bundles and the Hecke graph}

\begin{definition}
Define $\mathfrak{P}_C$ to be the subcomplex of the building
$\mathfrak{X}_C$ consisting of cells whose stabilizer in
$\op{GL}_3(k[C])$ contains a non-central subgroup. This
is called the \emph{parabolic subcomplex}, and its
$\op{GL}_3(k[C])$-equivariant homology
$$\widehat{\op{H}}_\bullet(\op{GL}_3(k[C]),M):=
\op{H}_\bullet^{\op{GL}_3(k[C])}(\mathfrak{P}_C,M)$$ 
is called \emph{parabolic homology}.
\end{definition}

By a result of Atiyah, a vector bundle $E$ on $\overline{C}$ is stable
if and only if $\op{Aut}(E)\cong k^\times$ consists of
scalars. Therefore, a $0$-simplex is in $\op{GL}_3(k[C])\backslash
\mathfrak{P}_C$ if and only if it corresponds to a bundle which is not
stable, i.e., is different from $E(3,1)$ and $E(3,2)$.

The complex $\op{GL}_3(k[C])\backslash\mathfrak{X}_C$ can be described
in terms of certain moduli spaces of vector bundles. For this, let me
define  a graph called ``Hecke graph'', the name chosen because the
graph is related to Hecke operators connecting moduli
spaces of vector bundles.

\begin{definition}
Let $C$ be an affine elliptic curve in Weierstrass normal form. 
Define the coloured \emph{Hecke graph} $\Gamma(C,k)$ to be 
$$
\mathcal{M}_{2,1}(\overline{C})(k)\stackrel{\alpha}{\longleftarrow}
\mathcal{M}_{2,0}(\overline{C})(k)\stackrel{\beta}{\longrightarrow}
\mathcal{M}_{3,\mathcal{O}}(\overline{C})(k).
$$ 
The map $\alpha:\op{Sym}^2(\overline{C})\to\op{Jac}(\overline{C})$ is the
determinant, and the map $\beta:\op{Sym}^2(\overline{C})\to\mathbb{P}^2$ is
the symmetric square of the (given choice of)  covering
$\overline{C}\to\mathbb{P}^1$ branched at four points. 
\end{definition}

The vertices in
$\mathcal{M}_{2,1}(\overline{C})\cong\op{Jac}(\overline{C})$
correspond to rank three 
bundles of the form
$E=E_{\mathcal{L}}(2,1)\oplus\mathcal{L}^{-1}$, and we have
$\op{Aut}(E)\cong k^\times\times k^\times$. 
The vertices in $\mathcal{M}_{3,\mathcal{O}}(\overline{C})\cong\mathbb{P}^2$
correspond to rank three bundles over $\overline{C}$ which
geometrically split as direct sums of three line bundles of degree $0$. 
There are five types of such bundles, depending on the structure of
the scheme-theoretic fiber of $\beta$ over $x\in\mathbb{P}^2$:
\begin{itemize}
\item Split bundles
  $E=\mathcal{L}_1\oplus\mathcal{L}_2\oplus\mathcal{L}_3$ with pairwise
  non-isomorphic summands, in which case $\op{Aut}(E)\cong
  (k^\times)^3$.
\item Split bundles
  $E=\mathcal{L}_1\oplus\mathcal{L}_1\oplus\mathcal{L}_2$ with
  $\mathcal{L}_1\not\cong \mathcal{L}_2$, in which case
  $\op{Aut}(E)\cong \op{GL}_2(k)\times k^\times$.
\item Split bundles
  $E=\mathcal{L}\oplus\mathcal{L}\oplus\mathcal{L}$ in which case
  $\op{Aut}(E)\cong\op{GL}_3(k)$. 
\item Partially split bundles
  $E=\mathcal{A}\oplus\det\mathcal{A}^{-1}$ with $\mathcal{A}$ an
  indecomposable geometrically split bundle corresponding to a
  degree 2 point $a\in\op{Jac}(\overline{C})$ with residue field
  $k(a)$. In this case, $\op{Aut}(E)\cong k(a)^\times\times
  k^\times$. 
\item Indecomposable but geometrically split bundles $E=\mathcal{B}$
  corresponding to a degree 3 point $b\in\op{Jac}(\overline{C})$ with
  residue field $k(b)$. In this case, $\op{Aut}(E)\cong k(b)^\times$. 
\end{itemize}

The set
$\mathcal{M}_{2,0}(\overline{C})\cong\op{Sym}^2(\overline{C})$ of
edges corresponds to bundles of 
the form $E=E_{\mathcal{L}}(2,0)\oplus\mathcal{L}^{-1}$. Up to a
unipotent subgroup, $\op{Aut}(E)\cong k^\times\times k^\times$.

\begin{theorem}
\label{thm:parquotient}
\begin{enumerate}
\item 
After a suitable subdivision of $\Gamma(C,k)$, the above assignment
gives rise to a simplicial weak equivalence 
$$\mu:\Gamma(C,k)\to
\op{GL}_3(k[C])\backslash\mathfrak{P}_C.$$ 
\item The assignment $\op{Aut}:\Gamma(C,k)\mapsto \mathcal{G}r$ turns
  $\Gamma(C,k)$ into a graph of groups. If $k$ is infinite or
  $\op{char}k$ is invertible in the coefficients, then the map $\mu$
  induces isomorphisms in Borel-equivariant homology of the
  developments of the respective complexes of groups.
\end{enumerate}
\end{theorem}

The proof of the theorem is a long list of tedious linear algebra
computations:  Atiyah's classification provides a list of all
bundles $E$ on $\overline{C}$ which are not stable, and for each of
these, one has to determine the action of $\op{Aut}(E)$ on
$\op{Lk}_{\mathfrak{X}_C}(e)$, where $e$ is a lift of the point of
$\op{GL}_3(k[C])\backslash\mathfrak{X}_C$ corresponding to
$E$. Eventually, this is the study of the action of subgroups of
$\op{GL}_3$ on projective homogeneous varieties $\op{GL}_3/P$ -
tedious, but nothing conceptually deep. A byproduct of the
computations are explicit representing cocycles for all rank three
bundles on $\overline{C}$ with trivial restriction to $C$. 

\begin{theorem}
\label{thm:parhlgy}
Denote by $\mathcal{E}=\mathcal{M}_{2,0}(\overline{C})(k)$ the set of
edges of $\Gamma(C,k)$, and by
$\mathcal{V}=\mathcal{M}_{2,1}(\overline{C})(k)
\cup\mathcal{M}_{3,\mathcal{O}}(\overline{C})(k)$
the set of vertices of $\Gamma(C,k)$.
There is a long exact sequence computing parabolic homology: 
$$
\cdots\to
\bigoplus_{E\in\mathcal{E}}\op{Aut}(E)
\stackrel{\alpha_\ast+\beta_\ast}{\longrightarrow} 
\bigoplus_{V\in\mathcal{V}}\op{Aut}(V)\to
\widehat{\op{H}}_\bullet(\op{GL}_3(k[C]),\mathbb{Z})\to\cdots
$$
where $\alpha_\ast$ and $\beta_\ast$ denote the induced maps on
homology for the inclusions of edge groups into vertex groups.
\end{theorem}

This follows from Theorem~\ref{thm:parquotient} by applying the
spectral sequence for $\op{GL}_3(k[C])$-equivariant homology of
$\mathcal{P}_C$. The identification with the Hecke graph of groups
shows that the spectral sequence reduces to a long exact sequence
exhibiting parabolic homology as cone of the inclusion of edge groups
into vertex groups. 

\begin{corollary}
Replacing $\op{Aut}(E)$ by $\op{Aut}(E)/k^\times$ in
Theorem~\ref{thm:parhlgy} provides a long exact sequence computing
$\op{H}_\bullet(\op{PGL}_3(k[C]),\mathbb{Z})$ for $\bullet\geq 3$.
\end{corollary}

This follows because the building for $\op{GL}_3(k[C])$ has dimension
$2$, and all the cells in $\op{PGL}_3(k[C])\backslash\mathfrak{X}_C$
which are not in the parabolic subcomplex have trivial stabilizer.

It is then fairly straightforward to compute Hilbert-Poincar{\'e}
series for the homology
$\op{H}_\bullet(\op{PGL}_3(\mathbb{F}_q[C]),\mathbb{F}_\ell)$ with
$\ell\nmid q$.

\section{Stable bundles and the Steinberg module}

It remains to describe the part of
$\op{GL}_3(k[C])\backslash\mathfrak{X}_C$ which is not contained in
the parabolic subcomplex. As noted earlier, the only $0$-simplices not
in $\op{GL}_3(k[C])\backslash\mathfrak{P}_C$ are the points $x_i$,
$i=1,2$ corresponding to the stable bundles $E(3,i)$. Denote by 
$\mathfrak{S}_C=\op{Star}_{\op{GL}_3(k[C])\backslash \mathfrak{X}_C}(x_1)\cup
\op{Star}_{\op{GL}_3(k[C])\backslash \mathfrak{X}_C}(x_2)$ the union
of the stars of $x_i$. Note that the boundary
$\op{Lk}_{\mathfrak{X}_C}(x_i)$ of these stars can be identified as
spherical building for $\op{GL}_3(k)$. 

\begin{theorem}
\label{thm:cusp}
The complex $\mathfrak{S}_C$  is a double cone over the 
spherical building for $\op{GL}_3(k)$. The inclusion 
$\mathfrak{S}_C\subseteq\op{GL}_3(k[C])\backslash\mathfrak{X}_C$ is a
weak equivalence of simplicial sets, hence
$\op{GL}_3(k[C])\backslash\mathfrak{X}_C$ has the homotopy type of a
suspension of the spherical building for $\op{GL}_3(k)$.
\end{theorem}

\begin{theorem}
\label{thm:hlgy}
\begin{enumerate}
\item
The inclusion of the parabolic subcomplex gives rise to an exact
sequence of complexes of $\op{GL}_3(k[C])$-modules 
$$
0\to\op{C}_\bullet(\mathfrak{P}_C)\to
\op{C}_\bullet(\mathfrak{X}_C)\to
\op{C}_\bullet(\mathfrak{X}_C/\mathfrak{P}_C)\to 0.
$$
\item
The quotient
$\mathfrak{Q}_C=\op{GL}_3(k[C])\backslash\mathfrak{X}_C/\mathfrak{P}_C$ 
has the homotopy type of a wedge $\mathfrak{S}_C\vee \Sigma\Gamma(C,k)$, and 
\begin{eqnarray*}
&&\op{H}_\bullet^{\op{GL}_3(k[C])}(\mathfrak{X}_C/\mathfrak{P}_C,\mathbb{Z})\\
&\cong &
\op{H}_\bullet(\mathfrak{Q}_C,\mathbb{Z})\otimes
\op{H}_\bullet(k^\times,\mathbb{Z})\\
&\cong&
\op{H}_\bullet(k^\times,\mathbb{Z})\oplus
\left(\op{H}_{\bullet-2}(k^\times,\mathbb{Z})\otimes
  \left(\op{St}_3(k)\oplus \op{H}_1(\Gamma(C,k),\mathbb{Z})\right)\right)
\end{eqnarray*}
where $\op{St}_3(k)$ is the Steinberg module for $\op{GL}_3(k)$. 
\item The  long exact homology sequence associated to the sequence in
  1) takes the form 
$$
\cdots\to
\widehat{\op{H}}_\bullet(\op{GL}_3(k[C]),\mathbb{Z})\to 
\op{H}_\bullet(\op{GL}_3(k[C]),\mathbb{Z})\to 
$$
$$
\to
\op{H}_\bullet(k^\times,\mathbb{Z})\oplus
\left(\op{H}_{\bullet-2}(k^\times,\mathbb{Z})\otimes
  \left(\op{St}_3(k)\oplus \op{H}_1(\Gamma(C,k),\mathbb{Z})\right)\right)
\to\cdots
$$
\end{enumerate}
\end{theorem}

The theorem follows directly from the description of the quotient in
Theorems~\ref{thm:parquotient} and \ref{thm:cusp}. The map
$\op{H}_\bullet(\op{GL}_3(k[C]),\mathbb{Z})\to\op{H}_\bullet(k^\times,\mathbb{Z})$
is the one induced from the determinant. Moreover, it is easy to see
that the restriction of the boundary map to
$$
\left(\op{H}_{\bullet-2}(k^\times,\mathbb{Z})\otimes
\op{H}_1(\Gamma(C,k),\mathbb{Z})\right)
\to \widehat{\op{H}}_{\bullet-1}(\op{GL}_3(k[C]),\mathbb{Z})
$$
is injective. However, at the moment I do not have a formula for
the differential 
$$
\left(\op{H}_{\bullet-2}(k^\times,\mathbb{Z})\otimes
  \op{St}_3(k)\right)\to
\widehat{\op{H}}_{\bullet-1}(\op{GL}_3(k[C]),\mathbb{Z}).
$$

For $\op{PGL}_3(k[C])$, there is a similar formula. It is in fact
easier, as the $\op{PGL}_3(k[C])$-equivariant homology of the quotient 
$\mathfrak{X}_C/\mathfrak{P}_C$ reduces to
$\op{St}_3(k)\oplus\op{H}_1(\Gamma(C,k),\mathbb{Z})$ in degree $2$ and
is trivial in degrees $\geq 3$.

\begin{corollary}
Let $k=\mathbb{F}_q$. Then
$\op{dim}_{\mathbb{Q}}\op{H}_2(\op{GL}_3(k[C]),\mathbb{Q})=q^3$. 
\end{corollary}

By the results of Harder \cite{harder}, we know that $\op{H}_0$ and
$\op{H}_2$ are the only possibly non-trivial rational homology groups
in this case. However, the explicit dimension computation seems to be
new. The above result implies the existence of a huge number of
``non-detectable'' cohomology classes in
$\op{H}^2(\op{GL}_3(k[C]),\mathbb{Z})$, i.e., classes which restrict
trivially to the diagonal matrices or any finite subgroups of
$\op{GL}_3(k[C])$.  This is yet another case where the function field
analogue of Quillen's conjecture on cohomology of $S$-arithmetic
groups fails.

\section{Algebraic K-theory of elliptic curves}

Of course, this research was started in the hope that it might have
applications to algebraic K-theory of elliptic curves. Although at
present, there is no explicit K-theoretic consequence of the above
computations, looking at the building may provide an alternative
approach to the construction of the motivic weight two complex for
elliptic curves: part (3) of Theorem~\ref{thm:hlgy} gives rise to a
map
$$
\sigma:\op{St}_3(k)\to\left(\widehat{\op{H}}_1(\op{GL}_3(k[C]),\mathbb{Z})/
\op{H}_1(\Gamma,\mathbb{Z})\right). 
$$

The cokernel of $\sigma$ is identified with
$\ker\left(\det:\op{H}_1(\op{GL}_3(k[C]))\to
  \op{H}_1(k^\times)\right)$, and stabilization results for the
homology of linear groups tell us that $\op{H}_1(\op{GL}_3(k[C]))\cong
\op{H}_1(\op{GL}_\infty(k[C]))$, and hence
$$
\ker\left(\det:\op{H}_1(\op{GL}_3(k[C]))\to
  \op{H}_1(k^\times)\right)
\cong \op{SK}_1(k[C]).
$$
In particular, an explicit computation of $\sigma$ (which at present
is not yet available) would provide an alternative presentation of
$\op{SK}_1(k[C])$.  This might shed light on Vaserstein's conjecture
which predicts $\op{SK}_1(k[C])$ to be torsion if $k$ is a number
field (or equivalently, predicts
$\op{SK}_1(\overline{\mathbb{Q}}[C])=0$). 

The kernel of $\sigma$ is identified with
$\op{H}_2(\op{GL}_3(k[C]))/\widehat{\op{H}}_2(\op{GL}_3(k[C]))$, and
stabilization results tell us that the natural map
$\op{H}_2(\op{GL}_3(k[C]))\to
\op{H}_2(\op{GL}_\infty(k[C]))\cong\op{K}_2(k[C])$ is
surjective. Again, this provides an alternative way of constructing
elements in $\op{K}_2(k[C])$. 

Putting these things together, we see that there is a (yet
unspecified) subgroup $\mathcal{R}_C\subseteq\ker\sigma\subseteq
\op{St}_3(k)$ (very likely depending on the curve $C$) such that the 
following sequence is exact
$$
0\to \op{K}^{\{2\}}_2(k[C])\to 
\op{St}_3(k)/\mathcal{R}_C
\stackrel{\sigma}{\longrightarrow}
\left(\widehat{\op{H}}_1(\op{GL}_3(k[C]),\mathbb{Z})/ 
\op{H}_1(\Gamma,\mathbb{Z})\right)\to
\op{SK}_1(k[C])\to 0
$$
where 
$$
\op{K}^{\{2\}}_2(k[C])=\op{H}_2(\op{GL}_\infty(k[E]))/
\op{Im}(\widehat{\op{H}}_2(\op{GL}_3(k[C])))
$$
denotes the ``rank two'' quotient of $\op{K}_2(k[C])$. By
stabilization, $\op{H}_2(\op{GL}_4(k[C]))\cong\op{K}_2(k[C])$, hence
the set $\mathcal{R}$ of relations could be determined by similar
investigations for the action of $\op{GL}_4(k[C])$ on the associated
building. 

I would expect a strong relation between the complex $[\sigma]$ and
the motivic weight two complex of \cite[Theorem
1.5]{goncharov:levin}. Recall that the theorem of Goncharov and Levin
states that the following is (after tensoring with $\mathbb{Q}$) an
exact sequence
$$
0\to \op{Tor}(k^\times,\op{Jac}(\overline{C}))\to
\op{H}^0(\overline{C},\mathcal{K}_2)/\op{K}_2(k)\to$$
$$
\to B_3^\ast(\overline{C})\to 
k^\times\otimes\op{Jac}(\overline{C})\to
\ker(\op{H}^1(\overline{C},\mathcal{K}_2)\to k^\times)\to 0
$$
so that $\tau:B_3^\ast(\overline{C})_{\mathbb{Q}}\to
k^\times\otimes\op{Jac}(\overline{C})_{\mathbb{Q}}$ is a motivic
weight two complex for the elliptic curve $\overline{C}$. See
\cite{goncharov:levin} for notation and definitions of elliptic Bloch
groups. 

There is an obvious conjecture relating the complex of
Goncharov-Levin (which is related to the complete curve
$\overline{C}$) to the complex $[\sigma]$ described in this note
(which is related to the open curve $C$):

\begin{conjecture}
I expect that there is a morphism of complexes
\begin{center}
  \begin{minipage}[c]{10cm}
    \xymatrix{
      B_3^\ast(\overline{C})_{\mathbb{Q}}\ar[r]\ar[d]_\tau &
      \left(\op{St}_3(k)/\mathcal{R}_C\right)_{\mathbb{Q}} \ar[d]^\sigma \\
      \left(k^\times\otimes\op{Jac}(\overline{C})\right)_{\mathbb{Q}}\ar[r]
      & 
      \left(\widehat{\op{H}}_1(\op{GL}_3(k[C]),\mathbb{Q})/ 
        \op{H}_1(\Gamma,\mathbb{Q})\right)      
    }
  \end{minipage}
\end{center}
which is natural in $k$ and $C$, such that the induced maps on
homology 
$$
\op{K}_2^{[2]}(\overline{C})_{\mathbb{Q}}\to
\op{K}_2^{\{2\}}(C)_{\mathbb{Q}}\quad\textrm{ and } 
\quad
\op{K}_1^{[2]}(\overline{C})_{\mathbb{Q}}\to \op{SK}_1(C)_{\mathbb{Q}}
$$
 agree with the natural restriction maps associated
to $C\hookrightarrow\overline{C}$. 
\end{conjecture}

If true, the conjecture would imply that the complex $[\sigma]$ is an
integral refinement of the motivic weight two complex. 

For the map
$k^\times\otimes\op{Jac}(\overline{C})\to
\left(\widehat{\op{H}}_1(\op{GL}_3(k[C]),\mathbb{Q})/
\op{H}_1(\Gamma,\mathbb{Q})\right)$, an element $u\otimes
\mathcal{L}$  naturally determines an element
$u\in\op{Aut}(\mathcal{L}\oplus\mathcal{L}^{-1})^{\op{ab}}$ (scaling
by $u$ on the first summand), and the image under the composition 
$\op{Aut}(\mathcal{L}\oplus\mathcal{L}^{-1})^{\op{ab}}\hookrightarrow
\op{H}_1(\op{GL}_2(k[C]))\to\op{H}_1(\op{GL}_3(k[C]))$ lies in
$\widehat{\op{H}}_1(\op{GL}_3(k[C]))$ and is independent of all choices.
For the map
$B_3^\ast(\overline{C})_{\mathbb{Q}}\to(\op{St}_3(k)/\mathcal{R}_C)_{\mathbb{Q}}$,
a divisor $D\in \mathbb{Z}[\overline{C}(k)]$ should be mapped to some
linear combination of $2$-simplices corresponding to elementary
transformations on bundles $\mathcal{V}\oplus\det\mathcal{V}^{-1}$
with $\det\mathcal{V}^{-1}$ concentrated in $D$. 
The crucial thing missing before the above can be made precise is the
explicit description of the map $\sigma$ and the relations
$\mathcal{R}_C$. 

As a final remark, one should note the structural similarity between
the above construction and the motivic weight two
complex for fields (due to Bloch-Wigner, Suslin, Dupont-Sah, etc):
$$
0\to \widetilde{\op{Tor}}(\mu_k,\mu_k)\to \op{K}_3^{\op{ind}}(k)\to \mathcal{P}(k)\to\Lambda^2(k^\times)\to
\op{K}_2(k)\to 0. 
$$
The complex itself arises from a differential in the spectral sequence
computing the $\op{GL}_2(k)$-equivariant homology of the \v Cech
resolution of $\mathbb{P}^1(k)$; the part $\Lambda^2(k^\times)$ arises
from the homology of stabilizer subgroups, and the scissors congruence
part $\mathcal{P}(k)$ is a quotient of
$\mathbb{Z}[\mathbb{P}^1]\cong\op{St}_2(k)$ modulo the
dilogarithmic functional equations. All this looks very similar to
what's done for elliptic curves in the present note.

The motivic weight two complexes (both for fields as well as for
elliptic curves) exhibit a relation between buildings for
$\op{GL}_n$ and cycle descriptions for algebraic K-theory which is not
yet sufficiently understood.

\end{document}